\documentclass[12pt]{article}
\usepackage{amssymb, amsmath, url, graphicx, setspace, geometry}
\usepackage{calrsfs}
\usepackage{wasysym}
\usepackage{xcolor}
\usepackage{tikz}

\def\3{\subset }
\def\4{\subseteq }
\def\<{\left<}
\def\>{\right>}

\def\bit{\begin{itemize}}
\def\eit{\end{itemize}}
\def\3{\subset }
\def\4{\subseteq }

\def\0{\leqno}

\def\barr{\begin{array}}
\def\earr{\end{array}}

\def\Z{{\rlap{$\kern2pt{\rm Z}$}{\rm Z}\,}}

\title{\bf A graph related to the sum of element orders of a finite group}
\author{Mihai-Silviu Lazorec}
\date{August 6, 2021}

\begin{document}
\maketitle

\begin{abstract}
A finite group is called $\psi$-divisible iff $\psi(H)|\psi(G)$ for any subgroup $H$ of a finite group $G$. Here, $\psi(G)$ is the sum of element orders of $G$. For now, the only known examples of such groups are the  cyclic ones of square-free order. The existence of non-abelian $\psi$-divisible groups still constitutes an open question. The aim of this paper is to make a connection between the $\psi$-divisibility property and graph theory. Hence, for a finite group $G$, we introduce a simple undirected graph called the $\psi$-divisibility graph of $G$. We denote it by $\psi_G$. Its vertices are the non-trivial subgroups of $G$, while two distinct vertices $H$ and $K$ are adjacent iff $H\subset K$ and $\psi(H)|\psi(K)$ or $K\subset H$ and $\psi(K)|\psi(H)$. We prove that $G$ is $\psi$-divisible iff $\psi_G$ has a universal (dominating) vertex. Also, we study various properties of $\psi_G$, when $G$ is a finite cyclic group. The choice of restricting our study to this specific class of groups is motivated in the paper. 
\end{abstract}



\section{Introduction}

In what follows given a simple undirected graph $\mathcal{G}$ we denote its vertex set by $V(\mathcal{G})$ and its edge set by $E(\mathcal{G})$. For an integer $n\geq 2$, let $C_n$ and $\tau(n)$ be the finite cyclic group of order $n$ and the number of positive divisors of $n$, respectively. For a finite group $G$, we denote by $\pi(G)$ and $L(G)$ the set containing the prime divisors of $|G|$ and the subgroup lattice of $G$, respectively, while the order of an element $x\in G$ is denoted by $o(x)$.

Even though one can go back to the 19th century to check Cayley's work \cite{20}, the connection between group and graph theories became more popular after 1950. There are a lot of graphs defined especially on a finite group $G$ and, without entering into much detail, we recall the following ones:
\begin{itemize}
\item \textit{the commuting graph} $Com(G)$ of $G$ in which $V(Com(G))=G$ and, for two distinct vertices $x$ and $y$, we have $xy\in E(Com(G))$ iff $xy=yx$ (see \cite{12, 17}); we note that this graph is often redefined such that $V(Com(G))=G\setminus Z(G)$ (see \cite{2, 13, 28, 32, 38});
\item \textit{the power graph} $P(G)$ of $G$ in which $V(P(G))=G$ and, for two distinct vertices $x$ and $y$, $xy\in E(P(G))$ iff $x\in\langle y\rangle$ or $y\in\langle x\rangle$ (see \cite{6, 14, 16, 18, 26, 33}); if one removes the vertex associated with the identity of $G$, then the obtained graph is called \textit{the proper power graph} of $G$ (see \cite{24, 39});
\item \textit{the enhanced power graph} $EP(G)$ of $G$ in which $V(P(G))=G$ and, for two distinct vertices $x$ and $y$, $xy\in E(EP(G))$ iff $\langle x,y\rangle$ is a cyclic subgroup (see \cite{1, 10, 48}); if one removes the vertex corresponding to the trivial element of $G$, then the obtained graph is called \textit{the deleted enhanced power graph} (see \cite{11}) or \textit{the cyclic graph} of $G$ (see \cite{22, 23});
\item \textit{the Gruenberg-Kegel} graph $GK(G)$ of $G$ in which $V(GK(G))=\pi(G)$ and, for two distinct vertices $p$ and $q$, $pq\in E(GK(G))$ iff $G$ contains an element $x$ such that $o(x)=pq$ (see \cite{19, 35, 47}); this graph is also called \textit{the prime graph} of $G$;
\end{itemize}
Paper \cite{15} is a recent survey containing some remarkable results concerning the graphs above. 

Another idea to obtain graphs related to a finite group $G$ is to use different subsets of $L(G)$ as being the vertex set and define the adjacency relation in a specific way. Some examples would be the intersection graph of $G$ (see \cite{27, 37, 40, 42, 44}), the join graph of $G$ (see \cite{3, 36}) or the recent  factorization graph of $G$ (see \cite{25}).

During the last decade, the sum of element orders of a finite group $G$ has proven to be a tool that may be used to characterize the nature and structure of $G$ (see \cite{4, 5, 7, 8, 9, 30, 31, 45, 46} for further details). We denote the sum of element orders of $G$ by $\psi(G)$. Its definition is straightforward as its name suggests, i.e.:
$$\psi(G)=\sum\limits_{x\in G}o(x).$$
In \cite{29}, the concept of $\psi$-divisible group was introduced. More exactly, we say that a finite group is $\psi-divisible$ if and only if $\psi(H)|\psi(G)$ for all $H\in L(G)$. The main result of the same paper states that a finite abelian group $G$ is $\psi$-divisible if and only if $G$ is cyclic of square-free order. The connection between $\psi$-divisibility and the square-free order property of finite groups was further investigated in \cite{34}. Among others, the same paper highlights the following result which refers to any finite group:\\

\textbf{Theorem 1.1.} \textit{Let $G$ be a finite group. Then all subgroups of $G$ are $\psi$-divisible if and only if $G$ is cyclic of square-free order.}\\

In this paper we introduce a graph which is related to the concepts outlined above. So, for a finite group $G$, we denote by $\psi_G$ a graph called the $\psi$-divisibility graph of $G$. Its vertex set is $V(\psi_G)=L(G)\setminus\lbrace \lbrace 1\rbrace\rbrace$, where $\lbrace 1\rbrace$ is the trivial subgroup of $G$. For two distinct vertices $H$ and $K$ we have $HK\in E(\psi_G)$ iff $H\subset K$ and $\psi(H)|\psi(K)$ or $K\subset H$ and $\psi(K)|\psi(H)$. We remove the vertex corresponding to $\lbrace 1\rbrace$ as it would be a universal (or dominating) vertex in our graph. As we previously saw, removing some potential vertices is a common idea in graph theory. In our case, by removing the trivial subgroup, some properties of the $\psi$-divisibility graph are more interesting to study (for instance, the connectivity). To avoid the generation of a $\psi$-divisibility graph isomorphic to the null graph, we work only with non-trivial finite groups even though we always omit to write ``non-trivial". 

We mention that the existence of $\psi$-divisible groups beyond the class of finite abelian groups is still an open problem. It would be ideal if one could say if and how the $\psi$-divisibility graph may be used to approach this problem. In this regard, we prove that a finite group $G$ is $\psi$-divisible if and only if its $\psi$-divisibility graph $\psi_G$ has a universal vertex. Another trouble is that it is quite difficult to generate the $\psi$-divisibility graph for any finite group. Taking into account this fact and Theorem 1.1, which provides a connection between the $\psi$-divisibility property and a part of the finite cyclic groups, most of our results concern the study of some properties of the $\psi$-divisibility graph associated especially with a finite cyclic group. More exactly, we are interested in connectivity (in this case we also discuss about diameter and number of connected components), vertex degrees, cycles, bipartiteness, girth, trees and planarity. Some open questions are outlined throughout the paper.

\section{Main results}

Let $G$ be a finite group. Since the vertices of the $\psi$-divisibility graph $\psi_G$ of $G$ are the non-trivial subgroups of $G$, it is obvious that the divisor graph of $|G|$ is a subgraph of $\psi_G$. The first result of this Section highlights some results which help us with drawing the edges of a $\psi$-divisibility graph. The main tools would be items $i)$ and $iii)$ of the following lemma. For the ease of writing, when we work with cyclic groups, we write $\psi(n)$ instead of $\psi(C_n)$.\\  

\textbf{Lemma 2.1.} \textit{Let $G_1, G_2$ be finite groups, $p$ be a prime number and let $a, b$ be positive integers. Then the following hold:
\begin{itemize}
\item[i)] $\psi(G_1\times G_2)=\psi(G_1)\psi(G_2)$ if and only if $(|G_1|, |G_2|)=1$ (i.e. $\psi$ is multiplicative);
\item[ii)] $\psi(p^a)=\frac{p^{2a+1}+1}{p+1}$;
\item[iii)] $\psi(p^a)|\psi(p^b)\Longleftrightarrow 2a+1|2b+1$;
\item[iv)] $(\psi(p), \psi(p^2))=(\psi(p), \psi(p^3))=(\psi(p^2), \psi(p^3))=1$.
\end{itemize}}

\textbf{Proof.} For items \textit{i)} and \textit{ii)}, one can check Lemma 2.1 of \cite{5} and Lemma 2.9 \textit{(1)} of \cite{30}, respectively. 

For item \textit{iii)}, we use some fundamental properties of the integer divisibility. We have
\begin{align*}
\psi(p^a)|\psi(p^b)&\Longleftrightarrow p^{2a+1}+1|p^{2b+1}+1\Longleftrightarrow -p^{2a+1}-1|-p^{2b+1}-1 \\ & \Longleftrightarrow (-p)^{2a+1}-1|(-p)^{2b+1}-1\Longleftrightarrow 2a+1|2b+1.
\end{align*}

Finally, concerning \textit{iv)}, we show that $(\psi(p), \psi(p^3))=1$. Let $d=(\psi(p), \psi(p^3))$. We have $d|(p^4-p)\psi(p)$ and $d|\psi(p^3)$. Since $(p^4-p)\psi(p)=\psi(p^3)-1$, we get $d|1$ and the conclusion follows. Using similar ideas, one can show that $(\psi(p), \psi(p^2))=(\psi(p^2), \psi(p^3))=1$.
\hfill\rule{1,5mm}{1,5mm}\\

In what concerns the completeness of the $\psi$-divisibility graph of an arbitrary finite group, it is easy to prove the following result:\\

\textbf{Proposition 2.2.} \textit{Let $G$ be a finite group. Then $\psi_G$ is complete if and only if $G\cong C_p$, where $p$ is a prime number.}

\textbf{Proof.} Let $G$ be a finite group and let $H\cong C_p$ be a subgroup of $G$, where $p\in\pi(G)$. Assume that $\psi_G$ is a complete graph. If $G$ contains a subgroup $K$ or order $p^2$, then $K\cong C_{p^2}$ or $K\cong C_p^2$. Since $H$ and $K$ must be adjacent, it follows that $\psi(p)|\psi(p^2)$ or $\psi(p)|\psi(C_p^2)$. The first case leads to a contradiction due to Lemma 2.1, \textit{iv)}. In the second case, it is easy to show that $\psi(p)|p$, again a contradiction. It follows that $G$ is of square-free order. Suppose that $|\pi(G)|\geq 2$. Then $L\cong C_q$ is another subgroup of $G$, where $q\in \pi(G)$. Due to the definition of the adjacency relation, we cannot draw the edge $HL$ of $\psi_G$, a contradiction. It follows that $G$ is of square-free order and $\pi(G)=\lbrace p\rbrace$, so $G\cong C_p$.

Conversely, we have $\psi_G\cong K_1$ and this fact completes our proof.
\hfill\rule{1,5mm}{1,5mm}\\

It is a difficult task to obtain an explicit result concerning the vertex degrees of the $\psi$-divisibility graph of an arbitrary finite group. We can say something if we restrict our study to the class of finite cyclic groups. \\

\textbf{Proposition 2.3.} \textit{Let $G$ be a finite cyclic group of order $n$.
\begin{itemize}
\item[i)] If $n=p^\alpha$, where $p\in\pi(G)$ and $\alpha\geq 1$, then $$d(C_{p^\beta})=\tau(2\beta+1)+\bigg\lfloor\frac{1}{2}\cdot\bigg(\frac{2\alpha+1}{2\beta+1}-1\bigg)\bigg\rfloor-2, \ \forall \ \beta\in\lbrace 1,2,\ldots,\alpha\rbrace;$$
\item[ii)] If $n=p_1^{\alpha_1}p_2^{\alpha_2}\ldots p_k^{\alpha_k}$, where $k\geq 2, p_i\in\pi(G)$ and $\alpha_i\geq 1$, then 
$$d(C_{p_i^{\alpha_i}})=\prod\limits_{j=1,j\ne i}^k(\alpha_j+1)+\tau(2\alpha_i+1)-3, \ \forall \ i\in\lbrace 1, 2,\ldots, k\rbrace.$$
\end{itemize}}
\textbf{Proof.} \textit{i)} Let $G\cong C_{p^{\alpha}}$, where $p$ is a prime and $\alpha\geq 1$, and let $\beta, \gamma\in\lbrace 1, 2,\ldots, \alpha\rbrace$, with $\beta\ne\gamma$. According to Lemma 2.1, \textit{iii)}, $C_{p^{\beta}}$ is adjacent with $C_{p^{\gamma}}$ if and only if $2\gamma+1>1$ is a proper divisor of $2\beta+1$ or $2\gamma+1$ is a multiple of $2\beta+1$ in the range $(2\beta+1, 2\alpha+1]$. As we know, the number of such multiples is $\lfloor\frac{2\alpha+1}{2\beta+1}\rfloor-1$. But, we need to take into account only the odd multiples, so we actually get the quantity  $\lfloor\frac{1}{2}\cdot(\frac{2\alpha+1}{2\beta+1}-1)\rfloor$. Putting all things together, we obtain
$$d(C_{p^\beta})=\tau(2\beta+1)+\bigg\lfloor\frac{1}{2}\cdot\bigg(\frac{2\alpha+1}{2\beta+1}-1\bigg)\bigg\rfloor-2.$$

\textit{ii)} Let $G\cong C_n$, where $n=p_1^{\alpha_1}p_2^{\alpha_2}\ldots p_k^{\alpha_k}$, $k\geq 2, p_i\in\pi(G)$ and $\alpha_i\geq 1$, for all $i\in\lbrace 1,2,\ldots, k\rbrace$. Fix $i\in\lbrace 1,2,\ldots, k\rbrace$. Due to how we defined the adjacency relation for $\psi_G$, the vertex $C_{p_i^{\alpha_i}}$ may be adjacent only with vertices contained in the subset $A_1\cup A_2$ of $V(\psi_G)$, where
$$A_1=L(C_{p_i^{\alpha_i}})\setminus\lbrace\lbrace 1\rbrace\rbrace \text{ and } A_2=\lbrace H\in L(G) \ | \ C_{p_i^{\alpha_i}}\subsetneq H\rbrace.$$
As we previously explained, the vertex $C_{p_i^{\alpha_i}}$ is adjacent with any vertex $C_{p^{\gamma}}\in A_1$ if and only if $2\gamma+1>1$ is a proper divisor of $2\alpha_i+1$. Since $L(G)$ is a decomposable lattice, due to Lemma 2.1, \textit{i)}, we deduce that $C_{p_i^{\alpha_i}}$ is adjacent with any vertex in  $A_2$. Then 
$$d(C_{p_i^{\alpha_i}})=\tau(2\alpha_i+1)-2+|A_2|=\prod\limits_{j=1,j\ne i}^k(\alpha_j+1)+\tau(2\alpha_i+1)-3.$$
\hfill\rule{1,5mm}{1,5mm}\\

Note that item \textit{ii)} of the previous result concerns only the degrees of the vertices associated with the Sylow subgroups of $G$. It would be interesting to obtain some similar formulas for all other remaining vertices. 

As we explained in the first Section, most of the subsequent results concern only the class of finite cyclic groups. Let $p$ and $q$ be two distinct primes. For the ease of writing, we outline the following 3 sets of divisibility relations:
$$(P_1) \left\{\begin{array}{ll}
\psi(p)|\psi(q)  \\
\psi(p)|\psi(q^2) \\
\psi(q)|\psi(p)\\
\psi(q)|\psi(p^2)
\end{array} \right. (P_2) \left\{\begin{array}{ll}
\psi(p)|\psi(q^2) \text{ \ and \ } \psi(q)|\psi(p^2)  \\
\psi(p)|\psi(q^2) \text{ \ and \ } \psi(q)|\psi(p^3) \\
\psi(p^2)|\psi(q^2) \text{ \ and \ } \psi(q)|\psi(p^3)
\end{array} \right. $$
$$(P_3) \left\{\begin{array}{ll}
\psi(p)|\psi(q^2) \text{ \ and \ } \psi(q)|\psi(p^2)  \\
\psi(p)|\psi(q^3) \text{ \ and \ } \psi(q)|\psi(p^3) \\
\psi(p)|\psi(q^3) \text{ \ and \ } \psi(q)|\psi(p^2) \\
\psi(p)|\psi(q^2) \text{ \ and \ } \psi(q)|\psi(p^3) \\
\psi(p^2)|\psi(q^3) \text{ \ and \ } \psi(q^2)|\psi(p^3) \\
\psi(p^2)|\psi(q^3) \text{ \ and \ } \psi(q)|\psi(p^3) \\
\psi(p)|\psi(q^3) \text{ \ and \ } \psi(q^2)|\psi(p^3) \\
\psi(p)|\psi(q^3) \text{ \ and \ } \psi(q^2)|\psi(p^2) \\
\psi(p^2)|\psi(q^2) \text{ \ and \ } \psi(q)|\psi(p^3)
\end{array} \right.$$
Before continuing we say a few words on these relations. Let $P$ be the set of the first $10^4$ primes. Concerning $(P_1)$, by using SageMath \cite{41}, we see that we can generate a lot of examples of pairs $(p,q)\in P\times P$ such that $\psi(p)|\psi(q)$. More exactly, there are 11631 such examples: $(2, 5), (2, 11), \ldots, (3181, 24841)$. The condition $\psi(p)|\psi(q^2)$ is more restrictive since $(p, q)\in\lbrace (151, 33469), (181, 14407), (181, 44483),$ $(181, 66851), (571, 24203), (571, 62591)\rbrace.$ As expected, any of the conditions included into $(P_2)$ or $(P_3)$ is even more restrictive since two divisibility relations must hold at the same time. Our code returns an example only for $\psi(p^2)|\psi(q^2) \text{ \ and \ } \psi(q)|\psi(p^3)$, this being $(p, q)=(2, 7)$. Obviously $(p, q)=(7, 2)$ would be an example such that $\psi(p)|\psi(q^3) \text{ \ and \ } \psi(q^2)|\psi(p^2)$. Hence, in general, it seems that there are high chances that none of the $(P_2)$ and, respectively, $(P_3)$ conditions holds.\\

One would say that cycle graphs in graph theory are equally important as cyclic groups in group theory. Hence, the next objective is to classify the finite cyclic groups whose $\psi$-divisibility graphs are cycles.\\

\textbf{Proposition 2.4.} \textit{Let $G$ be a finite cyclic group. Then $\psi_G$ is a cycle if an only if $G\cong C_{p^2q^2}$, where $p, q$ are prime numbers, and none of the ($P_1$) conditions holds.}

\textbf{Proof.} Let $G$ be a finite cyclic group and let $p\in \pi(G)$. Suppose that $\psi_G$ is a cycle. Then $\psi_G$ is a 2-regular graph. If $|\pi(G)|\geq 3$, then $d(C_p)\geq 3$, a contradiction. Hence $|\pi(G)|\in\lbrace 1,2\rbrace$. If $|\pi(G)|=1$, according to Lemma 2.1, \textit{iii)}, to draw only two edges incident with $C_p$, we need two odd multiples of 3. We deduce that $G\cong C_{p^\alpha}$, where $\alpha\in\lbrace 7, 8, 9\rbrace$. Indeed, in all these 3 cases, we can draw the edges $C_pC_{p^4}$ and $C_pC_{p^7}$. But, also in all these cases, one can easily check that $d(C_{p^5})=0$, a contradiction. 

Consequently $|\pi(G)|=2$, so there is an additional prime $q$ such that $G\cong C_{p^{\alpha}q^{\beta}}$, where $\alpha \geq 1, \beta\geq 1$. If $\alpha\geq 3$ or $\beta\geq 3$, we use Lemma 2.1, \textit{i)}, to obtain $d(C_q)\geq 3$ or $d(C_p)\geq 3$, a contradiction. If $\alpha=1$ or $\beta=1$, then $\psi_G$ has at least one end-vertex, again a contradiction. Hence $G\cong C_{p^2q^2}$. Finally, if at least one of the $(P_1)$ conditions are met, we get $d(C_p)\geq 3$ or $d(C_q)\geq 3$, a contradiction. 

Conversely, if $G\cong C_{p^2q^2}$ and all $(P_1)$ conditions do not hold, we can draw the $\psi$-divisibility graph below (Figure 1) and finish our proof. 
\hfill\rule{1,5mm}{1,5mm} 
\tikzset{every picture/.style={line width=0.75pt}} 
\begin{center}
\begin{tikzpicture}[x=0.75pt,y=0.75pt,yscale=-1,xscale=1]

\draw   (216.81,32.98) .. controls (216.81,25.81) and (223.53,20) .. (231.82,20) .. controls (240.11,20) and (246.83,25.81) .. (246.83,32.98) .. controls (246.83,40.14) and (240.11,45.95) .. (231.82,45.95) .. controls (223.53,45.95) and (216.81,40.14) .. (216.81,32.98) -- cycle ;
\draw   (397.56,54.06) .. controls (397.56,46.6) and (404.56,40.54) .. (413.19,40.54) .. controls (421.83,40.54) and (428.83,46.6) .. (428.83,54.06) .. controls (428.83,61.53) and (421.83,67.58) .. (413.19,67.58) .. controls (404.56,67.58) and (397.56,61.53) .. (397.56,54.06) -- cycle ;
\draw   (138,68.66) .. controls (138,61.19) and (145,55.14) .. (153.64,55.14) .. controls (162.27,55.14) and (169.27,61.19) .. (169.27,68.66) .. controls (169.27,76.12) and (162.27,82.17) .. (153.64,82.17) .. controls (145,82.17) and (138,76.12) .. (138,68.66) -- cycle ;
\draw   (310,33.52) .. controls (310,26.05) and (317,20) .. (325.63,20) .. controls (334.27,20) and (341.27,26.05) .. (341.27,33.52) .. controls (341.27,40.98) and (334.27,47.03) .. (325.63,47.03) .. controls (317,47.03) and (310,40.98) .. (310,33.52) -- cycle ;
\draw   (189.29,117.86) .. controls (189.29,110.39) and (196.29,104.34) .. (204.92,104.34) .. controls (213.56,104.34) and (220.56,110.39) .. (220.56,117.86) .. controls (220.56,125.32) and (213.56,131.37) .. (204.92,131.37) .. controls (196.29,131.37) and (189.29,125.32) .. (189.29,117.86) -- cycle ;
\draw   (381.92,132.46) .. controls (381.92,124.99) and (388.92,118.94) .. (397.56,118.94) .. controls (406.19,118.94) and (413.19,124.99) .. (413.19,132.46) .. controls (413.19,139.92) and (406.19,145.97) .. (397.56,145.97) .. controls (388.92,145.97) and (381.92,139.92) .. (381.92,132.46) -- cycle ;
\draw   (285.6,139.48) .. controls (285.6,132.02) and (292.6,125.97) .. (301.24,125.97) .. controls (309.87,125.97) and (316.88,132.02) .. (316.88,139.48) .. controls (316.88,146.95) and (309.87,153) .. (301.24,153) .. controls (292.6,153) and (285.6,146.95) .. (285.6,139.48) -- cycle ;
\draw   (463.23,98.93) .. controls (463.23,91.47) and (470.23,85.42) .. (478.86,85.42) .. controls (487.5,85.42) and (494.5,91.47) .. (494.5,98.93) .. controls (494.5,106.4) and (487.5,112.45) .. (478.86,112.45) .. controls (470.23,112.45) and (463.23,106.4) .. (463.23,98.93) -- cycle ;
\draw    (153.64,82.17) -- (193.66,108.13) ;
\draw    (220.56,117.86) -- (285.6,139.48) ;
\draw    (316.88,139.48) -- (381.92,132.46) ;
\draw    (153.64,55.14) -- (216.81,32.43) ;
\draw    (248.08,32.43) -- (310,33.52) ;
\draw    (341.27,33.52) -- (397.56,54.06) ;
\draw    (413.19,132.46) -- (478.86,112.45) ;
\draw    (424.45,63.79) -- (468.86,88.66) ;

\draw (142.48,60.89) node [anchor=north west][inner sep=0.75pt]  [font=\fontsize{0.73em}{0.88em}\selectfont] [align=left] {$\displaystyle C_{pq^2}$};
\draw (196.3,109.15) node [anchor=north west][inner sep=0.75pt]  [font=\fontsize{0.73em}{0.88em}\selectfont] [align=left] {$\displaystyle C_{p}$};
\draw (289.36,130.78) node [anchor=north west][inner sep=0.75pt]  [font=\fontsize{0.73em}{0.88em}\selectfont] [align=left] {$\displaystyle C_{pq}$};
\draw (388.93,124.83) node [anchor=north west][inner sep=0.75pt]  [font=\fontsize{0.73em}{0.88em}\selectfont] [align=left] {$\displaystyle C_{q}$};
\draw (465.61,90.77) node [anchor=north west][inner sep=0.75pt]  [font=\fontsize{0.73em}{0.88em}\selectfont] [align=left] {$\displaystyle C_{p^{2} q}$};
\draw (401.94,45.9) node [anchor=north west][inner sep=0.75pt]  [font=\fontsize{0.73em}{0.88em}\selectfont] [align=left] {$\displaystyle C_{p^{2}}$};
\draw (310.26,26.98) node [anchor=north west][inner sep=0.75pt]  [font=\fontsize{0.73em}{0.88em}\selectfont] [align=left] {$\displaystyle C_{p^{2} q^{2}}$};
\draw (221.82,24.81) node [anchor=north west][inner sep=0.75pt]  [font=\fontsize{0.73em}{0.88em}\selectfont] [align=left] {$\displaystyle C_{q^{2}}$};
\end{tikzpicture}
\end{center}
\begin{center}
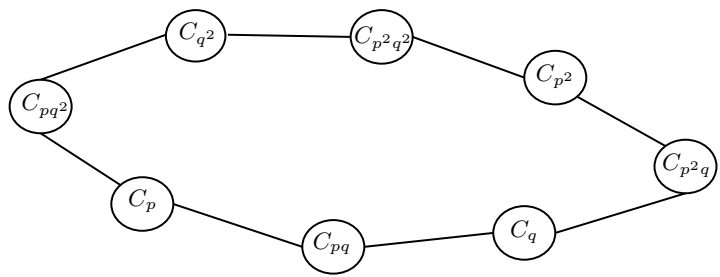
Figure 1. The graph $\psi_{C_{p^2q^2}}$ if none of the $(P_1)$ conditions holds
\end{center}
 
A well-known result in graph theory states that cycles can be used to characterize bipartite graphs. More exactly, a non-trivial graph is bipartite if and only if it does not contain odd cycles. In what follows, we use this result while classifying the finite cyclic groups whose $\psi$-divisibility graph is bipartite.\\

\textbf{Proposition 2.5.} \textit{Let $G$ be a finite cyclic group. Then $\psi_G$ is bipartite if and only if one of the following holds:
\begin{itemize}
\item[i)] $G\cong C_{p^{\alpha}}$, where $p$ is a prime and $\alpha\in\lbrace 2, 3,\ldots, 12\rbrace$;
\item[ii)] $G\cong C_{pq}, G\cong C_{p^2q}, G\cong C_{p^2q^2}$ or $G\cong C_{p^3q}$, where $p,q$ are primes;
\item[iii)] $G\cong C_{p^3q^2}$, where $p, q$ are primes, and none of the $(P_2)$ conditions holds;
\item[iv)] $G\cong C_{p^3q^3}$, where $p, q$ are primes, and none of the $(P_3)$ conditions holds.
\end{itemize}}
\textbf{Proof.} Let $G$ be a finite cyclic group and suppose that $\psi_G$ is a bipartite graph. If $|\pi(G)|\geq 3$, then we can choose $p, q, r\in\pi(G)$ and use the multiplicativity of $\psi$ to build the triangle $(C_p, C_{pq}, C_{pqr}, C_p)$, a contradiction. Consequently, $|\pi(G)|\in \lbrace 1, 2\rbrace$. 

\textit{i)} If $|\pi(G)|=1$, then $G\cong C_{p^{\alpha}}$, where $p\in\pi(G)$, while $\alpha\geq 2$ since $V(\psi_G)$ is partitioned into two non-empty partite sets, so we need at least two elements in $V(\psi_G)$. If $\alpha\geq 13$, then by Lemma 2.1, \textit{iii)}, we may consider the triangle $(C_p, C_{p^4}, C_{p^{13}}, C_p)$, a contradiction. Hence $\alpha\in\lbrace 2, 3,\ldots, 12\rbrace$.

The converse holds since $\psi_G$ would be a subgraph of the bipartite graph below (Figure 2).
\begin{center}
\tikzset{every picture/.style={line width=0.75pt}} 

\begin{tikzpicture}[x=0.75pt,y=0.75pt,yscale=-1,xscale=1]

\draw   (27.33,34.98) .. controls (27.33,27.81) and (33.84,22) .. (41.88,22) .. controls (49.91,22) and (56.43,27.81) .. (56.43,34.98) .. controls (56.43,42.15) and (49.91,47.96) .. (41.88,47.96) .. controls (33.84,47.96) and (27.33,42.15) .. (27.33,34.98) -- cycle ;
\draw   (82.61,34.98) .. controls (82.61,27.81) and (89.13,22) .. (97.16,22) .. controls (105.2,22) and (111.71,27.81) .. (111.71,34.98) .. controls (111.71,42.15) and (105.2,47.96) .. (97.16,47.96) .. controls (89.13,47.96) and (82.61,42.15) .. (82.61,34.98) -- cycle ;
\draw   (137.32,36.54) .. controls (137.32,29.37) and (143.83,23.56) .. (151.87,23.56) .. controls (159.9,23.56) and (166.42,29.37) .. (166.42,36.54) .. controls (166.42,43.7) and (159.9,49.51) .. (151.87,49.51) .. controls (143.83,49.51) and (137.32,43.7) .. (137.32,36.54) -- cycle ;
\draw   (194.93,36.02) .. controls (194.93,28.85) and (201.45,23.04) .. (209.48,23.04) .. controls (217.52,23.04) and (224.03,28.85) .. (224.03,36.02) .. controls (224.03,43.18) and (217.52,48.99) .. (209.48,48.99) .. controls (201.45,48.99) and (194.93,43.18) .. (194.93,36.02) -- cycle ;
\draw   (251.97,37.05) .. controls (251.97,29.89) and (258.48,24.08) .. (266.52,24.08) .. controls (274.55,24.08) and (281.07,29.89) .. (281.07,37.05) .. controls (281.07,44.22) and (274.55,50.03) .. (266.52,50.03) .. controls (258.48,50.03) and (251.97,44.22) .. (251.97,37.05) -- cycle ;
\draw   (315.4,38.09) .. controls (315.4,30.93) and (321.92,25.11) .. (329.95,25.11) .. controls (337.99,25.11) and (344.5,30.93) .. (344.5,38.09) .. controls (344.5,45.26) and (337.99,51.07) .. (329.95,51.07) .. controls (321.92,51.07) and (315.4,45.26) .. (315.4,38.09) -- cycle ;
\draw   (25,100.91) .. controls (25,93.74) and (31.51,87.93) .. (39.55,87.93) .. controls (47.58,87.93) and (54.1,93.74) .. (54.1,100.91) .. controls (54.1,108.07) and (47.58,113.89) .. (39.55,113.89) .. controls (31.51,113.89) and (25,108.07) .. (25,100.91) -- cycle ;
\draw   (80.29,100.91) .. controls (80.29,93.74) and (86.8,87.93) .. (94.84,87.93) .. controls (102.87,87.93) and (109.39,93.74) .. (109.39,100.91) .. controls (109.39,108.07) and (102.87,113.89) .. (94.84,113.89) .. controls (86.8,113.89) and (80.29,108.07) .. (80.29,100.91) -- cycle ;
\draw   (134.99,102.46) .. controls (134.99,95.3) and (141.51,89.49) .. (149.54,89.49) .. controls (157.58,89.49) and (164.09,95.3) .. (164.09,102.46) .. controls (164.09,109.63) and (157.58,115.44) .. (149.54,115.44) .. controls (141.51,115.44) and (134.99,109.63) .. (134.99,102.46) -- cycle ;
\draw   (192.61,101.95) .. controls (192.61,94.78) and (199.12,88.97) .. (207.16,88.97) .. controls (215.19,88.97) and (221.7,94.78) .. (221.7,101.95) .. controls (221.7,109.11) and (215.19,114.92) .. (207.16,114.92) .. controls (199.12,114.92) and (192.61,109.11) .. (192.61,101.95) -- cycle ;
\draw   (249.64,102.98) .. controls (249.64,95.82) and (256.15,90.01) .. (264.19,90.01) .. controls (272.22,90.01) and (278.74,95.82) .. (278.74,102.98) .. controls (278.74,110.15) and (272.22,115.96) .. (264.19,115.96) .. controls (256.15,115.96) and (249.64,110.15) .. (249.64,102.98) -- cycle ;
\draw   (313.07,104.02) .. controls (313.07,96.85) and (319.59,91.04) .. (327.62,91.04) .. controls (335.66,91.04) and (342.17,96.85) .. (342.17,104.02) .. controls (342.17,111.19) and (335.66,117) .. (327.62,117) .. controls (319.59,117) and (313.07,111.19) .. (313.07,104.02) -- cycle ;
\draw    (41.88,47.96) -- (39.55,87.93) ;
\draw    (53.23,43.28) -- (137.03,95.2) ;
\draw    (48.57,46.4) -- (86.98,90.52) ;
\draw    (151.87,49.51) -- (149.54,89.49) ;
\draw    (97.16,47.96) -- (94.84,87.93) ;
\draw    (109.97,40.69) -- (196.97,92.08) ;

\draw (35.21,28.98) node [anchor=north west][inner sep=0.75pt]  [font=\fontsize{0.73em}{0.88em}\selectfont] [align=left] {$\displaystyle C_{p}$};
\draw (89.08,29.61) node [anchor=north west][inner sep=0.75pt]  [font=\fontsize{0.73em}{0.88em}\selectfont] [align=left] {$\displaystyle C_{p^{2}}$};
\draw (143.2,29.9) node [anchor=north west][inner sep=0.75pt]  [font=\fontsize{0.73em}{0.88em}\selectfont] [align=left] {$\displaystyle C_{p^{3}}$};
\draw (200.01,30.76) node [anchor=north west][inner sep=0.75pt]  [font=\fontsize{0.73em}{0.88em}\selectfont] [align=left] {$\displaystyle C_{p^{5}}$};
\draw (257.74,30.65) node [anchor=north west][inner sep=0.75pt]  [font=\fontsize{0.73em}{0.88em}\selectfont] [align=left] {$\displaystyle C_{p^{6}}$};
\draw (321.17,32.21) node [anchor=north west][inner sep=0.75pt]  [font=\fontsize{0.73em}{0.88em}\selectfont] [align=left] {$\displaystyle C_{p^{8}}$};
\draw (30.88,94.16) node [anchor=north west][inner sep=0.75pt]  [font=\fontsize{0.73em}{0.88em}\selectfont] [align=left] {$\displaystyle C_{p^{4}}$};
\draw (86.06,94.27) node [anchor=north west][inner sep=0.75pt]  [font=\fontsize{0.73em}{0.88em}\selectfont] [align=left] {$\displaystyle C_{p^{7}}$};
\draw (139.38,95.95) node [anchor=north west][inner sep=0.75pt]  [font=\fontsize{0.73em}{0.88em}\selectfont] [align=left] {$\displaystyle C_{p^{10}}$};
\draw (196.88,95.83) node [anchor=north west][inner sep=0.75pt]  [font=\fontsize{0.73em}{0.88em}\selectfont] [align=left] {$\displaystyle C_{p^{12}}$};
\draw (255.18,97.73) node [anchor=north west][inner sep=0.75pt]  [font=\fontsize{0.73em}{0.88em}\selectfont] [align=left] {$\displaystyle C_{p^{9}}$};
\draw (316.66,98.02) node [anchor=north west][inner sep=0.75pt]  [font=\fontsize{0.73em}{0.88em}\selectfont] [align=left] {$\displaystyle C_{p^{11}}$};
\end{tikzpicture}
\end{center}
\begin{center}
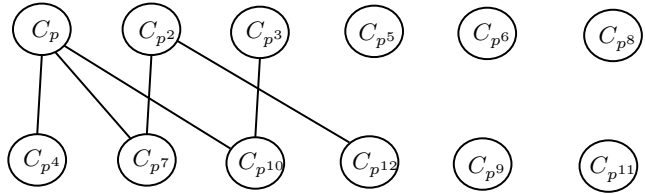
Figure 2. The graph $\psi_{C_{p^{12}}}$
\end{center}

If $|\pi(G)|=2$, then $G\cong C_{p^{\alpha}q^{\beta}}$, where $p, q\in\pi(G)$ and $\alpha\geq 1, \beta\geq 1$. If $\alpha\geq 4$ or $\beta\geq 4$, we can draw one of the triangles: $(C_p, C_{p^4}, C_{p^4q}, C_p)$, $(C_q, C_{q^4}, C_{pq^4}, C_{q})$, a contradiction. Hence $\alpha, \beta\in \lbrace 1,2,3\rbrace$. Due to symmetry, it suffices to check the cases $(\alpha, \beta)\in \lbrace (1,1), (2,1), (3,1), (2,2), (3,2), (3,3) \rbrace$.

\textit{ii)} If $(\alpha,\beta)\in\lbrace (1,1), (2,1), (2,2), (3,1)\rbrace$, then $G\cong C_{pq}, G\cong C_{p^2q}$, $G\cong C_{p^2q^2}$ or $G\cong C_{p^3q}$. The converse holds since $\psi_G$ would be a subgraph of one of the bipartite graphs below (Figure 3). We mention that some of the vertices of each partite set are not adjacent  due to Lemma 2.1, \textit{iv)}.
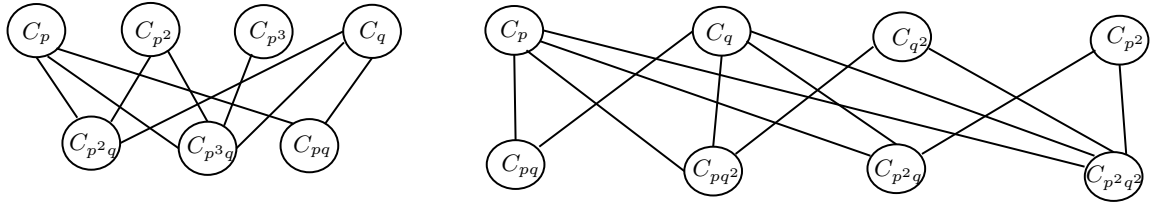
\begin{center}
\tikzset{every picture/.style={line width=0.75pt}} 

\begin{tikzpicture}[x=0.75pt,y=0.75pt,yscale=-1,xscale=1]

\draw   (200.61,27.42) .. controls (200.61,20.01) and (207.13,14) .. (215.16,14) .. controls (223.2,14) and (229.71,20.01) .. (229.71,27.42) .. controls (229.71,34.83) and (223.2,40.84) .. (215.16,40.84) .. controls (207.13,40.84) and (200.61,34.83) .. (200.61,27.42) -- cycle ;
\draw   (31.32,27.03) .. controls (31.32,19.62) and (37.83,13.61) .. (45.87,13.61) .. controls (53.9,13.61) and (60.42,19.62) .. (60.42,27.03) .. controls (60.42,34.44) and (53.9,40.45) .. (45.87,40.45) .. controls (37.83,40.45) and (31.32,34.44) .. (31.32,27.03) -- cycle ;
\draw   (88.93,26.49) .. controls (88.93,19.08) and (95.45,13.07) .. (103.48,13.07) .. controls (111.52,13.07) and (118.03,19.08) .. (118.03,26.49) .. controls (118.03,33.91) and (111.52,39.92) .. (103.48,39.92) .. controls (95.45,39.92) and (88.93,33.91) .. (88.93,26.49) -- cycle ;
\draw   (145.97,27.57) .. controls (145.97,20.16) and (152.48,14.15) .. (160.52,14.15) .. controls (168.55,14.15) and (175.07,20.16) .. (175.07,27.57) .. controls (175.07,34.98) and (168.55,40.99) .. (160.52,40.99) .. controls (152.48,40.99) and (145.97,34.98) .. (145.97,27.57) -- cycle ;
\draw   (169,85.12) .. controls (169,77.71) and (175.51,71.7) .. (183.55,71.7) .. controls (191.58,71.7) and (198.1,77.71) .. (198.1,85.12) .. controls (198.1,92.53) and (191.58,98.54) .. (183.55,98.54) .. controls (175.51,98.54) and (169,92.53) .. (169,85.12) -- cycle ;
\draw   (58.99,83.83) .. controls (58.99,76.42) and (65.51,70.41) .. (73.54,70.41) .. controls (81.58,70.41) and (88.09,76.42) .. (88.09,83.83) .. controls (88.09,91.25) and (81.58,97.25) .. (73.54,97.25) .. controls (65.51,97.25) and (58.99,91.25) .. (58.99,83.83) -- cycle ;
\draw   (117.64,86.44) .. controls (117.64,79.03) and (124.15,73.02) .. (132.19,73.02) .. controls (140.22,73.02) and (146.74,79.03) .. (146.74,86.44) .. controls (146.74,93.85) and (140.22,99.86) .. (132.19,99.86) .. controls (124.15,99.86) and (117.64,93.85) .. (117.64,86.44) -- cycle ;
\draw    (103.48,39.92) -- (83.5,73.01) ;
\draw    (154.5,38.89) -- (140.5,75.08) ;
\draw    (56,36) -- (175.5,74) ;
\draw    (215.16,40.84) -- (191.5,74) ;
\draw    (201,33) -- (146.74,86.44) ;
\draw    (200.61,27.42) -- (88.09,83.83) ;
\draw    (45.87,40.45) -- (66.5,71) ;
\draw    (51.5,40) -- (117.64,86.44) ;
\draw    (112,37) -- (132.19,73.02) ;
\draw   (577.61,31.87) .. controls (577.61,25.07) and (584.13,19.56) .. (592.16,19.56) .. controls (600.2,19.56) and (606.71,25.07) .. (606.71,31.87) .. controls (606.71,38.67) and (600.2,44.18) .. (592.16,44.18) .. controls (584.13,44.18) and (577.61,38.67) .. (577.61,31.87) -- cycle ;
\draw   (272.32,26.92) .. controls (272.32,20.12) and (278.83,14.61) .. (286.87,14.61) .. controls (294.9,14.61) and (301.42,20.12) .. (301.42,26.92) .. controls (301.42,33.73) and (294.9,39.24) .. (286.87,39.24) .. controls (278.83,39.24) and (272.32,33.73) .. (272.32,26.92) -- cycle ;
\draw   (376.93,27.35) .. controls (376.93,20.55) and (383.45,15.04) .. (391.48,15.04) .. controls (399.52,15.04) and (406.03,20.55) .. (406.03,27.35) .. controls (406.03,34.15) and (399.52,39.66) .. (391.48,39.66) .. controls (383.45,39.66) and (376.93,34.15) .. (376.93,27.35) -- cycle ;
\draw   (467.97,30.17) .. controls (467.97,23.37) and (474.48,17.86) .. (482.52,17.86) .. controls (490.55,17.86) and (497.07,23.37) .. (497.07,30.17) .. controls (497.07,36.97) and (490.55,42.48) .. (482.52,42.48) .. controls (474.48,42.48) and (467.97,36.97) .. (467.97,30.17) -- cycle ;
\draw   (465,96.74) .. controls (465,89.94) and (471.51,84.43) .. (479.55,84.43) .. controls (487.58,84.43) and (494.1,89.94) .. (494.1,96.74) .. controls (494.1,103.54) and (487.58,109.05) .. (479.55,109.05) .. controls (471.51,109.05) and (465,103.54) .. (465,96.74) -- cycle ;
\draw   (272.99,94.64) .. controls (272.99,87.84) and (279.51,82.33) .. (287.54,82.33) .. controls (295.58,82.33) and (302.09,87.84) .. (302.09,94.64) .. controls (302.09,101.44) and (295.58,106.96) .. (287.54,106.96) .. controls (279.51,106.96) and (272.99,101.44) .. (272.99,94.64) -- cycle ;
\draw   (372.64,97.95) .. controls (372.64,91.15) and (379.15,85.64) .. (387.19,85.64) .. controls (395.22,85.64) and (401.74,91.15) .. (401.74,97.95) .. controls (401.74,104.75) and (395.22,110.26) .. (387.19,110.26) .. controls (379.15,110.26) and (372.64,104.75) .. (372.64,97.95) -- cycle ;
\draw   (574.61,100.69) .. controls (574.61,93.88) and (581.13,88.37) .. (589.16,88.37) .. controls (597.2,88.37) and (603.71,93.88) .. (603.71,100.69) .. controls (603.71,107.49) and (597.2,113) .. (589.16,113) .. controls (581.13,113) and (574.61,107.49) .. (574.61,100.69) -- cycle ;
\draw    (286.87,39.24) -- (287.54,82.33) ;
\draw    (293,36.99) -- (372.64,97.95) ;
\draw    (300,32.4) -- (466.5,90.21) ;
\draw    (301.42,26.92) -- (574.5,95.71) ;
\draw    (376.93,27.35) -- (299.5,85.62) ;
\draw    (391.48,39.66) -- (387.19,85.64) ;
\draw    (404,32.4) -- (479.55,84.43) ;
\draw    (406.03,27.35) -- (579.5,90.21) ;
\draw    (496,36.07) -- (589.16,88.37) ;
\draw    (468,35.15) -- (398.5,89.29) ;
\draw    (592.16,44.18) -- (595.5,89.29) ;
\draw    (580,36.99) -- (492.5,89.29) ;

\draw (207.08,21.08) node [anchor=north west][inner sep=0.75pt]  [font=\fontsize{0.73em}{0.88em}\selectfont] [align=left] {$\displaystyle C_{q}$};
\draw (37.2,21.43) node [anchor=north west][inner sep=0.75pt]  [font=\fontsize{0.73em}{0.88em}\selectfont] [align=left] {$\displaystyle C_{p}$};
\draw (94.01,21.28) node [anchor=north west][inner sep=0.75pt]  [font=\fontsize{0.73em}{0.88em}\selectfont] [align=left] {$\displaystyle C_{p^{2}}$};
\draw (151.74,21.15) node [anchor=north west][inner sep=0.75pt]  [font=\fontsize{0.73em}{0.88em}\selectfont] [align=left] {$\displaystyle C_{p^{3}}$};
\draw (173.88,78.36) node [anchor=north west][inner sep=0.75pt]  [font=\fontsize{0.73em}{0.88em}\selectfont] [align=left] {$\displaystyle C_{pq}$};
\draw (120.18,79.14) node [anchor=north west][inner sep=0.75pt]  [font=\fontsize{0.73em}{0.88em}\selectfont] [align=left] {$\displaystyle C_{p^{3} q}$};
\draw (60.88,77.47) node [anchor=north west][inner sep=0.75pt]  [font=\fontsize{0.73em}{0.88em}\selectfont] [align=left] {$\displaystyle C_{p^{2} q}$};
\draw (584.08,23.72) node [anchor=north west][inner sep=0.75pt]  [font=\fontsize{0.73em}{0.88em}\selectfont] [align=left] {$\displaystyle C_{p^{2}}$};
\draw (278.2,20.33) node [anchor=north west][inner sep=0.75pt]  [font=\fontsize{0.73em}{0.88em}\selectfont] [align=left] {$\displaystyle C_{p}$};
\draw (382.01,22.03) node [anchor=north west][inner sep=0.75pt]  [font=\fontsize{0.73em}{0.88em}\selectfont] [align=left] {$\displaystyle C_{q}$};
\draw (475.74,24.7) node [anchor=north west][inner sep=0.75pt]  [font=\fontsize{0.73em}{0.88em}\selectfont] [align=left] {$\displaystyle C_{q^{2}}$};
\draw (466.88,90.92) node [anchor=north west][inner sep=0.75pt]  [font=\fontsize{0.73em}{0.88em}\selectfont] [align=left] {$\displaystyle C_{p^{2} q}$};
\draw (375.18,90.76) node [anchor=north west][inner sep=0.75pt]  [font=\fontsize{0.73em}{0.88em}\selectfont] [align=left] {$\displaystyle C_{pq^{2}}$};
\draw (278.88,88.26) node [anchor=north west][inner sep=0.75pt]  [font=\fontsize{0.73em}{0.88em}\selectfont] [align=left] {$\displaystyle C_{pq}$};
\draw (573.61,93.77) node [anchor=north west][inner sep=0.75pt]  [font=\fontsize{0.73em}{0.88em}\selectfont] [align=left] {$\displaystyle C_{p^{2} q^{2}}$};
\end{tikzpicture}
\end{center}
\begin{center}
Figure 3.
\end{center}

\textit{iv)} The two remaining choices for $(\alpha, \beta)$ may be investigated similarly, so we only treat the case $(\alpha,\beta)=(3,3)$ corresponding to the last item of our result. Therefore, $G\cong C_{p^3q^3}$. Since $\psi_G$ is bipartite, we are able to determine the partition $\lbrace V_1, V_2\rbrace$ of $V(\psi_G)$. Without loss of generality, we assume that $C_p\in V_1$. Then $C_{pq}, C_{pq^2}, C_{pq^3}$ must be placed in $V_2$ to avoid adjacency. By continuing this process, we get $V_1=\lbrace C_p, C_q, C_{p^2}, C_{q^2}, C_{p^3}, C_{q^3}\rbrace$ and $V_2=\lbrace C_{pq}, C_{p^2q^2}, C_{p^3q^3}, C_{pq^2}, C_{p^2q}, C_{pq^3}, C_{p^3q}, C_{p^2q^3},$ $C_{p^3q^2}\rbrace$. Note that, by construction, $V_1$ is an independent set of vertices. Since $V_2$ must satisfy the same property, we conclude that we can not draw any of the following potential edges: $C_{pq}C_{p^2q^2}, C_{pq}C_{p^3q^3}, C_{pq}C_{p^2q^3}, C_{pq}C_{p^3q^2}$, $C_{p^2q^2}C_{p^3q^3}$, $C_{p^2q}C_{p^3q^3}, C_{pq^2}C_{p^3q^3}$, $C_{pq^2}C_{p^2q^3},$ $C_{p^2q}C_{p^3q^2}$, i.e. none of the $(P_3)$ conditions holds. 

Conversely, if $G\cong C_{p^3q^3}$ and none of the $(P_3)$ conditions is met, then $\psi_G$ is isomorphic to a subgraph of the bipartite graph below (Figure 4). The edges marked with red mean that the vertices $C_p$ and $C_q$ may be adjacent with any of the vertices placed in $V_2$.
\hfill\rule{1,5mm}{1,5mm}
\begin{center}
\tikzset{every picture/.style={line width=0.75pt}} 

\begin{tikzpicture}[x=0.75pt,y=0.75pt,yscale=-1,xscale=1]

\draw   (86.67,109.24) .. controls (86.67,101.61) and (93.54,95.43) .. (102.02,95.43) .. controls (110.5,95.43) and (117.37,101.61) .. (117.37,109.24) .. controls (117.37,116.87) and (110.5,123.05) .. (102.02,123.05) .. controls (93.54,123.05) and (86.67,116.87) .. (86.67,109.24) -- cycle ;
\draw   (111.84,36.14) .. controls (111.84,28.51) and (118.71,22.33) .. (127.19,22.33) .. controls (135.67,22.33) and (142.54,28.51) .. (142.54,36.14) .. controls (142.54,43.77) and (135.67,49.96) .. (127.19,49.96) .. controls (118.71,49.96) and (111.84,43.77) .. (111.84,36.14) -- cycle ;
\draw   (196.85,34.37) .. controls (196.85,26.74) and (203.72,20.56) .. (212.2,20.56) .. controls (220.68,20.56) and (227.55,26.74) .. (227.55,34.37) .. controls (227.55,42) and (220.68,48.18) .. (212.2,48.18) .. controls (203.72,48.18) and (196.85,42) .. (196.85,34.37) -- cycle ;
\draw   (47.05,33.92) .. controls (47.05,26.3) and (53.92,20.11) .. (62.4,20.11) .. controls (70.88,20.11) and (77.75,26.3) .. (77.75,33.92) .. controls (77.75,41.55) and (70.88,47.74) .. (62.4,47.74) .. controls (53.92,47.74) and (47.05,41.55) .. (47.05,33.92) -- cycle ;
\draw   (21.11,109.35) .. controls (21.11,101.72) and (27.99,95.54) .. (36.46,95.54) .. controls (44.94,95.54) and (51.81,101.72) .. (51.81,109.35) .. controls (51.81,116.98) and (44.94,123.16) .. (36.46,123.16) .. controls (27.99,123.16) and (21.11,116.98) .. (21.11,109.35) -- cycle ;
\draw   (351.85,35.37) .. controls (351.85,27.74) and (358.72,21.56) .. (367.2,21.56) .. controls (375.68,21.56) and (382.55,27.74) .. (382.55,35.37) .. controls (382.55,43) and (375.68,49.18) .. (367.2,49.18) .. controls (358.72,49.18) and (351.85,43) .. (351.85,35.37) -- cycle ;
\draw   (272.85,33.37) .. controls (272.85,25.74) and (279.72,19.56) .. (288.2,19.56) .. controls (296.68,19.56) and (303.55,25.74) .. (303.55,33.37) .. controls (303.55,41) and (296.68,47.18) .. (288.2,47.18) .. controls (279.72,47.18) and (272.85,41) .. (272.85,33.37) -- cycle ;
\draw   (273.85,106.37) .. controls (273.85,98.74) and (280.72,92.56) .. (289.2,92.56) .. controls (297.68,92.56) and (304.55,98.74) .. (304.55,106.37) .. controls (304.55,114) and (297.68,120.18) .. (289.2,120.18) .. controls (280.72,120.18) and (273.85,114) .. (273.85,106.37) -- cycle ;
\draw   (204.85,108.09) .. controls (204.85,100.86) and (211.72,95) .. (220.2,95) .. controls (228.68,95) and (235.55,100.86) .. (235.55,108.09) .. controls (235.55,115.32) and (228.68,121.18) .. (220.2,121.18) .. controls (211.72,121.18) and (204.85,115.32) .. (204.85,108.09) -- cycle ;
\draw   (342.85,107.37) .. controls (342.85,99.74) and (349.72,93.56) .. (358.2,93.56) .. controls (366.68,93.56) and (373.55,99.74) .. (373.55,107.37) .. controls (373.55,115) and (366.68,121.18) .. (358.2,121.18) .. controls (349.72,121.18) and (342.85,115) .. (342.85,107.37) -- cycle ;
\draw   (142.85,107.37) .. controls (142.85,99.74) and (149.72,93.56) .. (158.2,93.56) .. controls (166.68,93.56) and (173.55,99.74) .. (173.55,107.37) .. controls (173.55,115) and (166.68,121.18) .. (158.2,121.18) .. controls (149.72,121.18) and (142.85,115) .. (142.85,107.37) -- cycle ;
\draw   (436.85,36.37) .. controls (436.85,28.74) and (443.72,22.56) .. (452.2,22.56) .. controls (460.68,22.56) and (467.55,28.74) .. (467.55,36.37) .. controls (467.55,44) and (460.68,50.18) .. (452.2,50.18) .. controls (443.72,50.18) and (436.85,44) .. (436.85,36.37) -- cycle ;
\draw   (417.85,106.37) .. controls (417.85,98.74) and (424.72,92.56) .. (433.2,92.56) .. controls (441.68,92.56) and (448.55,98.74) .. (448.55,106.37) .. controls (448.55,114) and (441.68,120.18) .. (433.2,120.18) .. controls (424.72,120.18) and (417.85,114) .. (417.85,106.37) -- cycle ;
\draw   (491.85,107.37) .. controls (491.85,99.74) and (498.72,93.56) .. (507.2,93.56) .. controls (515.68,93.56) and (522.55,99.74) .. (522.55,107.37) .. controls (522.55,115) and (515.68,121.18) .. (507.2,121.18) .. controls (498.72,121.18) and (491.85,115) .. (491.85,107.37) -- cycle ;
\draw   (547.85,107.37) .. controls (547.85,99.74) and (554.72,93.56) .. (563.2,93.56) .. controls (571.68,93.56) and (578.55,99.74) .. (578.55,107.37) .. controls (578.55,115) and (571.68,121.18) .. (563.2,121.18) .. controls (554.72,121.18) and (547.85,115) .. (547.85,107.37) -- cycle ;
\draw   (16.81,97.35) .. controls (16.81,92.93) and (20.4,89.35) .. (24.81,89.35) -- (578.5,89.35) .. controls (582.92,89.35) and (586.5,92.93) .. (586.5,97.35) -- (586.5,121.35) .. controls (586.5,125.77) and (582.92,129.35) .. (578.5,129.35) -- (24.81,129.35) .. controls (20.4,129.35) and (16.81,125.77) .. (16.81,121.35) -- cycle ;
\draw [color={rgb, 255:red, 208; green, 2; blue, 27 }  ,draw opacity=1 ]   (62.4,47.74) -- (61.5,90) ;
\draw    (198,40) -- (165.5,95) ;
\draw    (212.2,48.18) -- (220.2,95) ;
\draw    (219,45) -- (273.85,106.37) ;
\draw    (225,41) -- (342.85,107.37) ;
\draw    (227.55,34.37) -- (417.85,106.37) ;
\draw    (272.85,33.37) -- (113.5,102) ;
\draw    (279,43) -- (232.5,100) ;
\draw    (288.2,47.18) -- (296.5,95) ;
\draw    (295,44) -- (346.5,99) ;
\draw    (301,39) -- (423.5,97) ;
\draw    (351.85,35.37) -- (303.5,101) ;
\draw    (367.2,49.18) -- (433.2,92.56) ;
\draw    (382.55,35.37) -- (491.85,107.37) ;
\draw    (439,42) -- (373.55,107.37) ;
\draw    (452.2,50.18) -- (443.5,98) ;
\draw    (462,45) -- (552.5,98) ;
\draw [color={rgb, 255:red, 208; green, 2; blue, 27 }  ,draw opacity=1 ]   (127.19,49.96) -- (126.5,89) ;

\draw (52.84,26.33) node [anchor=north west][inner sep=0.75pt]  [font=\fontsize{0.73em}{0.88em}\selectfont] [align=left] {$\displaystyle C_{p}$};
\draw (117.84,30.03) node [anchor=north west][inner sep=0.75pt]  [font=\fontsize{0.73em}{0.88em}\selectfont] [align=left] {$\displaystyle C_{q}$};
\draw (201.09,27.7) node [anchor=north west][inner sep=0.75pt]  [font=\fontsize{0.73em}{0.88em}\selectfont] [align=left] {$\displaystyle C_{p^{2}}$};
\draw (143.3,100.92) node [anchor=north west][inner sep=0.75pt]  [font=\fontsize{0.73em}{0.88em}\selectfont] [align=left] {$\displaystyle C_{p^{2} q}$};
\draw (25.11,104.35) node [anchor=north west][inner sep=0.75pt]  [font=\fontsize{0.73em}{0.88em}\selectfont] [align=left] {$\displaystyle C_{pq}$};
\draw (275.3,100) node [anchor=north west][inner sep=0.75pt]  [font=\fontsize{0.73em}{0.88em}\selectfont] [align=left] {$\displaystyle C_{p^{3} q^{2}}$};
\draw (356.09,28.7) node [anchor=north west][inner sep=0.75pt]  [font=\fontsize{0.73em}{0.88em}\selectfont] [align=left] {$\displaystyle C_{p^{3}}$};
\draw (278.84,26) node [anchor=north west][inner sep=0.75pt]  [font=\fontsize{0.73em}{0.88em}\selectfont] [align=left] {$\displaystyle C_{q{^{2}}}$};
\draw (205.46,101.26) node [anchor=north west][inner sep=0.75pt]  [font=\fontsize{0.73em}{0.88em}\selectfont] [align=left] {$\displaystyle C_{p^{2} q{^{2}}}$};
\draw (343.3,100) node [anchor=north west][inner sep=0.75pt]  [font=\fontsize{0.73em}{0.88em}\selectfont] [align=left] {$\displaystyle C_{p^{2} q{^{3}}}$};
\draw (89.46,102.26) node [anchor=north west][inner sep=0.75pt]  [font=\fontsize{0.73em}{0.88em}\selectfont] [align=left] {$\displaystyle C_{pq{^{2}}}$};
\draw (440.09,31.7) node [anchor=north west][inner sep=0.75pt]  [font=\fontsize{0.73em}{0.88em}\selectfont] [align=left] {$\displaystyle C_{q^{3}}$};
\draw (418.3,101) node [anchor=north west][inner sep=0.75pt]  [font=\fontsize{0.73em}{0.88em}\selectfont] [align=left] {$\displaystyle C_{p^{3} q^{3}}$};
\draw (493.3,101) node [anchor=north west][inner sep=0.75pt]  [font=\fontsize{0.73em}{0.88em}\selectfont] [align=left] {$\displaystyle C_{p^{3} q}$};
\draw (551.3,102) node [anchor=north west][inner sep=0.75pt]  [font=\fontsize{0.73em}{0.88em}\selectfont] [align=left] {$\displaystyle C_{pq^{3}}$};
\end{tikzpicture}
\end{center}
\begin{center}
Figure 4.
\end{center}

Figure 2 clearly shows that, in general, the $\psi$-divisibility graph does not determine a finite group. For instance $\psi_{C_{2^{12}}}\cong \psi_{C_{3^{12}}}$, but $C_{2^{12}}\ncong C_{3^{12}}$. We continue by highlighting the possible values of $g(\psi_G)$, i.e. the values of the girth of the $\psi$-divisibility graph of a finite cyclic group $G$.\\

\textbf{Corollary 2.6.} \textit{Let $G$ be a finite cyclic group. Then $g(\psi_G)\in\lbrace 3,4,8,\infty\rbrace$.}

\textbf{Proof.} Let $G$ be a finite cyclic group. Using some of the reasoning outlined in the proof of Proposition 2.5, we can state that $g(\psi_G)=3$ in any of the following cases:
\begin{itemize}
\item[--] $|\pi(G)|\geq 3$;
\item[--] $G\cong C_{p^{\alpha}}$, where $p\in\pi(G)$ and $\alpha\geq 13$;
\item[--] $G\cong C_{p^{\alpha}q^{\beta}}$, where $p,q\in\pi(G)$ and $\alpha\geq 4$ or $\beta\geq 4$.
\end{itemize} 
If $G\cong C_{p^{\alpha}}$, where $p$ is a prime and $\alpha\leq 12$, then $\psi_G$ is a subgraph of the graph in Figure 2 and it follows that $g(\psi_G)=\infty$. If $G\cong C_{p^{\alpha}q^{\beta}}$, where $p,q\in \pi(G)$ and $1\leq\alpha,\beta\leq 3$, once again, due to symmetry it suffices to investigate the cases $(\alpha, \beta)\in \lbrace (1,1), (2,1), (3,1), (2,2), (3,2), (3,3) \rbrace$. 

Assume that $(\alpha, \beta)\in \lbrace (1,1), (2,1), (3,1)\rbrace$. Then $\psi_G$ contains no cycles, i.e. $g(\psi_G)=\infty$, excepting the following cases in which we have $g(\psi_G)=4$:
\begin{itemize}
\item[--] $(\alpha, \beta)=(2,1)$ and $\psi(p)|\psi(q)$; 
\item[--] $(\alpha, \beta)=(3,1)$ and $\psi(p)|\psi(q)$ or $\psi(p^2)|\psi(q)$.
\end{itemize}
Suppose that $(\alpha,\beta)=(2,2)$. As a consequence of Proposition 2.4, we obtain
\begin{align*}
g(\psi_G)=\begin{cases} 4, \ if \ at \ least \ one \ of \ the \ (P_1) \ conditions \ holds \\  8, \ if \ none \ of \ the \ (P_1) \ conditions \ holds  \end{cases}.
\end{align*}
Finally, assume that $(\alpha,\beta)\in \lbrace (3,2), (3,3)\rbrace$. It is clear that $3\leq g(\psi_G)\leq 8$ since $\psi_G$ contains a subgraph as the one in Figure 1. Since there are a lot of cases to handle due to which of the $(P_1), (P_2), (P_3)$ conditions hold, we used SageMath \cite{41} to check that our conclusion is true and finish our proof. 
\hfill\rule{1,5mm}{1,5mm}\\

The connectivity is one of the most relevant properties of a graph. When a new graph is introduced, one of the main questions is if it is connected or not. Further, based on the answer, one would be interested in determining the diameter or the number of connected components of the graph, respectively. Before proving our next result, let us consider the $\psi$-divisibility graph of $C_{p^n}$, where $n\geq 1$, and let $C_{p^{\alpha}}$, with $1\leq \alpha\leq n$ be one of its vertices. Then, according to Lemma 2.1,  \textit{iii)}, $C_{p^{\alpha}}$ is an isolated vertex of $\psi_{C_{p^n}}$ if and only if $2\alpha+1$ does not have a proper divisor $d>1$ or a multiple $m$ such that $2\alpha+1<m\leq 2n+1$. Hence, if we denote by $I(V(\psi_{C_{p^n}}))$ the set of isolated vertices of the $\psi$-divisibility graph $\psi_{C_{p^n}}$, then it is easy to determine its size. In what follows, we study the connectivity of the $\psi$-divisibility graph of a finite cyclic group $G$. There are some trivial cases such as $G\cong C_{p^n}$, where $p\in\pi(G)$ and $n\in \lbrace 1, 2, 3\rbrace$. If $n=1$, then $\psi_G$ is the trivial graph which is connected and its diameter is 0. If $n=2$ or $n=3$, then $\psi_G$ is disconnected and its number of connected components is 2 or 3, respectively.\\ 

\textbf{Theorem 2.7.} \textit{Let $G$ be a finite cyclic group.
\begin{itemize}
\item[i)] If $|\pi(G)|\geq 2$, then $\psi_G$ is connected. Moreover, $2\leq diam(\psi_G)\leq 4$;
\item[ii)] If $G\cong C_{p^n}$, where $p$ is a prime and $n\geq 4$, then $\psi_G$ is disconnected. In addition,  $k(\psi_G)=1+|I(V(\psi_G))|.$
\end{itemize}
}
\textbf{Proof.} \textit{i)} Let $G\cong C_{|G|}, |\pi(G)|=k\geq 2$ and let $|G|=p_1^{\alpha_1}p_2^{\alpha_2}\ldots p_k^{\alpha_k}$, where $p_l\in \pi(G)$ and $\alpha_l\geq 1$ for all $l\in\lbrace 1,2,\ldots, k\rbrace$. Take $i, j\in\lbrace 1,2,\ldots, k\rbrace$ and the divisors $x=p_1^{\beta_1}p_2^{\beta_2}\ldots p_k^{\beta_k}$ and $y=p_1^{\gamma_1}p_2^{\gamma_2}\ldots p_k^{\gamma_k}$ of $|G|$, where $\beta_l\leq\alpha_l, \gamma_l\leq \alpha_l$ for all $l\in\lbrace 1,2,\ldots, k\rbrace$ and $\beta_i>0, \gamma_j>0$. We use the multiplicativity of $\psi$ to draw a $C_x - C_y$ path. 

If $(x, y)=1$, then $(C_x, C_{xy}, C_y)$ is a $C_x - C_y$ path. Further, we assume that $(x, y)\ne 1$. If $C_x$ and $C_y$ are adjacent, we are done, so we also assume that we cannot draw the edge $C_xC_y$. If there is $l\in\lbrace 1,2,\ldots, k\rbrace$ such that at least one of $\beta_l$ and $\gamma_l$ is 0, then one can choose one of the following $C_x - C_y$ paths:
\begin{itemize}
\item[--] $(C_x, C_{xp_l}, C_{p_l}, C_{yp_l}, C_y)$, if $\beta_l=\gamma_l=0$;
\item[--] $(C_x, C_{xp_l^{\gamma_l}}, C_{p_l^{\gamma_l}}, C_y)$, if $\beta_l=0$ and $\gamma_l\ne 0$;
\item[--] $(C_x, C_{p_l^{\beta_l}}, C_{yp_l^{\beta_l}}, C_y)$, if $\beta_l\ne 0$ and $\gamma_l=0$.
\end{itemize}
If $\beta_l\ne 0$ and $\gamma_l\ne 0$ for all $l\in\lbrace 1,2,\ldots, k\rbrace$, then $(C_x, C_{p_i^{\beta_i}}, C_{p_1^{\gamma_1}\ldots p_{i-1}^{\gamma_{i-1}}p_i^{\beta_i}p_{i+1}^{\gamma_{i+1}}\ldots p_k^{\gamma_k}}$, $C_{p_1^{\gamma_1}\ldots p_{i-1}^{\gamma_{i-1}}p_{i+1}^{\gamma_{i+1}}\ldots p_k^{\gamma_k}}, C_y)$ is also a $C_x - C_y$ path. As a consequence of all these cases, we conclude that $\psi_G$ is a connected graph. 

Note that the length of all determined paths is at most 4. In addition,  $d(C_{p_1}, C_{p_2})=2$. Consequently, we have $2\leq diam(\psi_G)\leq 4$.   

\textit{ii)} Let $G\cong C_{p^n}$, where $p$ is a prime number and $n\geq 4$. To justify the disconnectedness of $\psi_G$, it is sufficient to show that $I(V(\psi_G))\ne\emptyset$. According to Bertrand's postulate, there is a prime $q=2z+1$ such that $n+1<q<2n+1$. Then $C_{p^z}$ is an isolated vertex in $\psi_G$. 

Further, in what concerns the number of connected components of $\psi_G$, it is sufficient to show that $k(\psi_G-I(V(\psi_G)))=1$. Hence, let $C_{p^{\alpha}}, C_{p^{\beta}}$ be two non-isolated vertices of $\psi_G$, where $\alpha, \beta\in \lbrace 1, 2,\ldots, n\rbrace$ and $\alpha\ne\beta$. Hence, we can choose two positive odd integers $s>1, t>1$ such that:
\begin{itemize}
\item[--] $2s+1$ is a proper divisor or a multiple of $2\alpha+1$, with $2\alpha+1<2s+1\leq 2n+1$ in the latter case;
\item[--] $2t+1$ is a proper divisor or a multiple of $2\beta+1$, with $2\beta+1<2t+1\leq 2n+1$ in the latter case.
\end{itemize}
To complete our proof, we must find a $C_{p^{\alpha}} - C_{p^{\beta}}$ path in our graph. Assume that $C_{p^{\alpha}}$ and $C_{p^{\beta}}$ are not adjacent and let $d_s=2u+1>1$ and $d_t=2v+1>1$ be the lowest divisors of $2s+1$ and $2t+1$, respectively. First, we prove that the vertex $C_{p^{\frac{d_sd_t-1}{2}}}$ exists in our graph. It suffices to show that $d_sd_t\leq 2n+1$. Without loss of generality, we assume that $d_s\leq d_t$. We have
$$d_sd_t\leq d_s\frac{2t+1}{d_t}=\frac{d_s}{d_t}(2t+1)\leq 2t+1\leq 2n+1.$$
Then, we are able to draw the $C_{p^{\alpha}} - C_{p^{\beta}}$ path $(C_{p^{\alpha}}, C_{p^s}, C_{p^u}, C_{p^{\frac{d_sd_t-1}{2}}}, C_{p^v}, C_{p^t}, C_{p^{\beta}})$, so $\psi_G-I(V(\psi_G))$ is a connected graph, as desired.
\hfill\rule{1,5mm}{1,5mm}\\

The disconnectedness of the $\psi$-divisibility graph of $C_{p^n}$, where $p$ is a prime and $n\geq 4$, can be justified by avoiding Bertrand's postulate. We insert a second proof of the result which uses a different idea: by assuming that the graph is connected, we are able to ``generate" more and more vertices. So, let $G\cong C_{p^n}$, where $p$ is a prime and $n\geq 4$. Assume that $\psi_G$ is a connected graph. As a first step, we observe that the vertices $C_{p_1}$ and $C_{p_1^2}$ are not adjacent but, since they are connected, it follows that $C_{p_1^7}\in V(\psi_G)$. This means that the graph contains the vertices $C_{p_1^{a_1}}$, where $a_1\in \lbrace 3, 4, 5, 6, 7\rbrace=A_1$. For the second step, we use the fact that $C_{p_1^6}$ and $C_{p_1^7}$ must be connected and this implies that $C_{p_1^{a_2}}\in V(\psi_G)$, where $a_2\in \lbrace 8, 9, \ldots, 97\rbrace=A_2$. At step $i\geq 3$, we choose the highest two elements of $A_{i-1}$, say $x_i$ and $y_i=x_i+1$. We observe that $X_i=2x_i+1$ and $Y_i=2y_i+1=X_i+2$ are relatively prime since $X_i$ is odd. Then $C_{p_1^{a_i}}\in V(\psi_G)$, where $a_i\in\lbrace y_i+1, y_i+2, \ldots, \frac{X_iY_i-1}{2}\rbrace=A_i$. The process continues by repeating the same ideas and it ``generates" a countably infinite set of vertices, a contradiction. Therefore, $\psi_G$ is a disconnected graph.

Note that the converses of items \textit{i)} and \textit{ii)} of Theorem 2.7 also hold. Using this remark, we are able to classify the finite cyclic groups whose $\psi$-divisibility graph is a tree.\\

\textbf{Corollary 2.8.} \textit{Let $G$ be a finite cyclic group. Then $\psi_G$ is a tree if and only if one of the following holds:
\begin{itemize}
\item[i)] $G\cong C_p$, where $p$ is a prime;
\item[ii)] $G\cong C_{pq}$, where $p, q$ are primes;
\item[iii)] $G\cong C_{p^2q}$, where $p, q$ are primes such that $\psi(p)\nmid \psi(q)$;
\item[iv)] $G\cong C_{p^3q}$, where $p, q$ are primes such that $\psi(p)\nmid \psi(q)$ and $\psi(p^2)\nmid \psi(q)$.
\end{itemize}}

\textbf{Proof.} Let $G$ be a finite cyclic group such that $\psi_G$ is a tree. It follows that $\psi_G$ is connected, so $G\cong C_p$, where $p$ is a prime, or $|\pi(G)|\geq 2$. In the latter case, according to the proof of Corollary 2.6, $\psi_G$ contains no cycles if and only if $G\cong C_{p^{\alpha}q^{\beta}}$, where $p,q\in\pi(G)$ and one of the following holds:
\begin{itemize}
\item[--] $(\alpha, \beta)=(1,1)$;
\item[--] $(\alpha, \beta)=(2,1)$ and $\psi(p)\nmid \psi(q)$;
\item[--] $(\alpha, \beta)=(3,1)$, $\psi(p)\nmid \psi(q)$ and $\psi(p^2)\nmid \psi(q)$.
\end{itemize}

The converse holds since $\psi_G$ is trivial for item \textit{i)}, the path $(C_p, C_{pq}, C_q)$ for item \textit{ii)}, the path $(C_p, C_{pq}, C_q, C_{p^2q}, C_{p^2})$ for item \textit{iii)} or the graph below (Figure 5) for item \textit{iv)}.
\hfill\rule{1,5mm}{1,5mm}
\begin{center}
\tikzset{every picture/.style={line width=0.75pt}} 

\begin{tikzpicture}[x=0.75pt,y=0.75pt,yscale=-1,xscale=1]

\draw   (261.67,68.24) .. controls (261.67,60.61) and (268.54,54.43) .. (277.02,54.43) .. controls (285.5,54.43) and (292.37,60.61) .. (292.37,68.24) .. controls (292.37,75.87) and (285.5,82.05) .. (277.02,82.05) .. controls (268.54,82.05) and (261.67,75.87) .. (261.67,68.24) -- cycle ;
\draw   (259.84,25.14) .. controls (259.84,17.51) and (266.71,11.33) .. (275.19,11.33) .. controls (283.67,11.33) and (290.54,17.51) .. (290.54,25.14) .. controls (290.54,32.77) and (283.67,38.96) .. (275.19,38.96) .. controls (266.71,38.96) and (259.84,32.77) .. (259.84,25.14) -- cycle ;
\draw   (260.85,108.37) .. controls (260.85,100.74) and (267.72,94.56) .. (276.2,94.56) .. controls (284.68,94.56) and (291.55,100.74) .. (291.55,108.37) .. controls (291.55,116) and (284.68,122.18) .. (276.2,122.18) .. controls (267.72,122.18) and (260.85,116) .. (260.85,108.37) -- cycle ;
\draw   (200.05,107.92) .. controls (200.05,100.3) and (206.92,94.11) .. (215.4,94.11) .. controls (223.88,94.11) and (230.75,100.3) .. (230.75,107.92) .. controls (230.75,115.55) and (223.88,121.74) .. (215.4,121.74) .. controls (206.92,121.74) and (200.05,115.55) .. (200.05,107.92) -- cycle ;
\draw   (201.11,67.35) .. controls (201.11,59.72) and (207.99,53.54) .. (216.46,53.54) .. controls (224.94,53.54) and (231.81,59.72) .. (231.81,67.35) .. controls (231.81,74.98) and (224.94,81.16) .. (216.46,81.16) .. controls (207.99,81.16) and (201.11,74.98) .. (201.11,67.35) -- cycle ;
\draw   (317.85,111.37) .. controls (317.85,103.74) and (324.72,97.56) .. (333.2,97.56) .. controls (341.68,97.56) and (348.55,103.74) .. (348.55,111.37) .. controls (348.55,119) and (341.68,125.18) .. (333.2,125.18) .. controls (324.72,125.18) and (317.85,119) .. (317.85,111.37) -- cycle ;
\draw   (315.85,70.37) .. controls (315.85,62.74) and (322.72,56.56) .. (331.2,56.56) .. controls (339.68,56.56) and (346.55,62.74) .. (346.55,70.37) .. controls (346.55,78) and (339.68,84.18) .. (331.2,84.18) .. controls (322.72,84.18) and (315.85,78) .. (315.85,70.37) -- cycle ;
\draw    (289.5,31) -- (325.5,58) ;
\draw    (275.19,38.96) -- (277.02,54.43) ;
\draw    (263,33) -- (222.5,55) ;
\draw    (216.46,81.16) -- (215.4,94.11) ;
\draw    (277.02,82.05) -- (276.2,94.56) ;
\draw    (331.2,84.18) -- (333.2,97.56) ;

\draw (208.84,99.33) node [anchor=north west][inner sep=0.75pt]  [font=\fontsize{0.73em}{0.88em}\selectfont] [align=left] {$\displaystyle C_{p}$};
\draw (265.84,17.03) node [anchor=north west][inner sep=0.75pt]  [font=\fontsize{0.73em}{0.88em}\selectfont] [align=left] {$\displaystyle C_{q}$};
\draw (267.09,99.7) node [anchor=north west][inner sep=0.75pt]  [font=\fontsize{0.73em}{0.88em}\selectfont] [align=left] {$\displaystyle C_{p^{2}}$};
\draw (263.3,61.92) node [anchor=north west][inner sep=0.75pt]  [font=\fontsize{0.73em}{0.88em}\selectfont] [align=left] {$\displaystyle C_{p^{2} q}$};
\draw (205.46,59.26) node [anchor=north west][inner sep=0.75pt]  [font=\fontsize{0.73em}{0.88em}\selectfont] [align=left] {$\displaystyle C_{pq}$};
\draw (318.3,63) node [anchor=north west][inner sep=0.75pt]  [font=\fontsize{0.73em}{0.88em}\selectfont] [align=left] {$\displaystyle C_{p^{3} q}$};
\draw (323.09,103.7) node [anchor=north west][inner sep=0.75pt]  [font=\fontsize{0.73em}{0.88em}\selectfont] [align=left] {$\displaystyle C_{p^{3}}$};
\end{tikzpicture}
\end{center}
\begin{center}
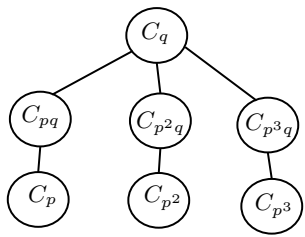
Figure 5. The graph $\psi_{C_{p^3q}}$ if $\psi(p)\nmid \psi(q)$ and $\psi(p^2)\nmid \psi(q)$
\end{center}

We saw that the diameter of a $\psi$-divisibility graph associated with a finite cyclic group takes low values: 2, 3 or 4. Without success, we tried to determine the conditions in which this graph is of diameter 2.  We insert an open problem concerning this aspect.\\

\textbf{Open problem.} \textit{Let $G$ be a finite cyclic group. Prove that if $diam(\psi_G)=2$, then $|\pi(G)|\geq 2$ and $G$ is of square-free order.}\\

The converse of the above statement clearly holds. The $\psi$-divisibility graph in Figure 5 is of diameter 4, while examples of $\psi$-divisibility graphs of diameter 3 are the ones associated with groups such as $C_{p^2q}$, where $p,q\in\pi(G)$ and $\psi(p)|\psi(q)$.\\

It is known that a graph $\mathcal{G}$ of order $|V(\mathcal{G})|\geq 3$ and size $|E(\mathcal{G})|$ is non-planar if $|E(\mathcal{G})|>3|V(\mathcal{G})|-6$. We use this result to show that most of the $\psi$-divisibility graphs of finite cyclic groups are non-planar.\\

\textbf{Proposition 2.9.} \textit{Let $G$ be a finite cyclic group such that $|\pi(G)|\geq 4$. Then $\psi_{G}$ is a non-planar graph.}

\textbf{Proof.} Let $G$ be a finite cyclic group with $|\pi(G)|\geq 4$. Then $\psi_{C_{pqrs}}$ is a subgraph of $\psi_{G}$, where $p, q, r, s\in\pi(G)$. It suffices to show that this subgraph is non-planar. Note that $|V(\psi_{C_{pqrs}})|=15$ and the vertex degrees with respect to $\psi_{C_{pqrs}}$ are:
$$d(C_p)=d(C_q)=d(C_r)=d(C_s)=d(C_{pqr})=d(C_{prs})=d(C_{pqs})=d(C_{qrs})=7;$$
$$d(C_{pq})=d(C_{pr})=d(C_{ps})=d(C_{qr})=d(C_{qs})=d(C_{rs})=5; d(C_{pqrs})=14.$$
Then  
$$|E(\psi_{C_{pqrs}})|=\frac{1}{2}\sum\limits_{v\in V(\psi_{C_{pqrs}})}d(v)=50>39=3|V(\psi_{C_{pqrs}})|-6,$$
so $\psi_{C_{pqrs}}$ is non-planar, as desired.
\hfill\rule{1,5mm}{1,5mm}\\

For sure there are a lot of problems related to graph theory that were not discussed in our paper. Some ideas for further research are outlined in the following paragraphs.\\

\textbf{Open problem.} Study other properties (regularity, clique number, automorphism group, spectrum, $L$-spectrum, number of spanning trees, chromatic number, vertex-connectivity, eulerianity, hamiltonicity, finite group recognition, etc.) of the $\psi$-divisibility graph of a finite cyclic group.\\

A starting point concerning the regularity would be Proposition 2.3. For a finite cyclic group $G$ such that  $\pi(G)=\lbrace p_1, p_2,\ldots, p_k\rbrace$, with $k\geq 1$, it is clear that $\omega(\psi_G)\geq |\pi(G)|$, where $\omega(\psi_G)$ is the clique number of $\psi_G$. This happens since the subgraph induced by the vertices $C_{p_1}, C_{p_1p_2}, \ldots, C_{p_1p_2\ldots p_k}$ is isomorphic to $K_{|\pi(G)|}$. 

Also, it would be interesting to study the $\psi$-divisibility graph beyond the class of finite cyclic groups. As starting points, we suggest to work with finite abelian $p$-groups of small rank, finite $p$-groups possessing a cyclic maximal subgroup (see Theorem 4.1, \cite{43}, vol. II) or finite groups having ``many" cyclic subgroups.\\

\textbf{Open problem.} Study the properties of the $\psi$-divisibility graph for specific classes of finite groups.\\

We end our paper with a result which establishes a connection between the $\psi$-divisibility property and the $\psi$-divisibility graph of an arbitrary finite group. In this way, one is able to identify $\psi$-divisible groups using graph theory.\\

\textbf{Theorem 2.10.} \textit{Let $G$ be a finite group. Then $G$ is $\psi$-divisible if and only if $\psi_G$ has a universal vertex.}

\textbf{Proof.} Let $G$ be a finite group. Obviously, if $G$ is $\psi$-divisible, then $G$ is a universal vertex of $\psi_G$. 

Conversely, assume that $\psi_G$ has a universal vertex. Due to how we defined the adjacency relation for our graph, it follows that the universal vertex is $G$ or it is a vertex corresponding to a breaking point in the subgroup lattice of $G$. In the first case, we are done since $G$ would be $\psi$-divisible. In the second case, $G$ would be isomorphic to one of the so called finite BP-groups (see \cite{21}). According to Theorem 1.1 of \cite{21} it follows that $G\cong C_{p^n}$, with $p\in\pi(G)$ and $n\geq 2$, or $G\cong Q_{2^n}$, where
$$Q_{2^n}=\langle x,y \ | \ x^{2^{n-1}}=y^4=1, yxy^{-1}=x^{2^{n-1}-1}\rangle, n\geq 3,$$
is the generalized quaternion group. 

If we assume that $G\cong C_{p^n}$, then we arrive at a contradiction since $\psi_G$ would be disconnected (check Theorem 2.7, \textit{ii)}, and the remarks preceding it concerning the cases $n\in\lbrace 2, 3\rbrace$). Suppose that $G\cong Q_{2^n}$. Then the universal vertex corresponds to the breaking point $H=\langle x^{2^{n-2}}\rangle\cong C_2$ of the subgroup lattice $L(Q_{2^n})$. Consider the vertex $K=\langle x^{2^{n-2}}\rangle\cong C_4$. Since $H$ is a universal vertex and $H\subset K$, it follows that $\psi(2)|\psi(4)$, so $3|11$, a contradiction. Thus the proof is complete.  
\hfill\rule{1,5mm}{1,5mm}\\

Obviously, we can use other graph theory concepts to restate Theorem 2.10. For instance we can say that a finite group $G$ is $\psi$-divisible if and only if $\gamma(\psi_G)=1$, where $\gamma(\psi_G)$ is the domination number of $\psi_G$. Recall that a finite abelian groups is $\psi$-divisible iff it is cyclic of square-free order. Then, as a consequence of Theorem 2.10, one can state the following result:\\

\textbf{Corollary 2.11.} \textit{Let $G$ be a finite abelian group. Then the following conditions are equivalent:
\begin{itemize}
\item[i)] $G$ is $\psi$-divisible;
\item[ii)] $G$ is cyclic of square-free order;
\item[iii)] $\psi_G$ has a universal vertex.
\end{itemize}}

\bigskip\noindent {\bf Acknowledgements.} This work was supported by a grant of the "Alexandru Ioan Cuza" University of Iasi, within the Research Grants program, Grant UAIC, code GI-UAIC-2021-01.

\vspace*{3ex}
\small
\hfill
\begin{minipage}[t]{6cm}
Mihai-Silviu Lazorec \\
Faculty of  Mathematics \\
"Al.I. Cuza" University \\
Ia\c si, Romania \\
e-mail: {\tt silviu.lazorec@uaic.ro}
\end{minipage}

\end{document}